\documentclass[a4paper,12pt]{article}
\usepackage{amsmath, amssymb, amsthm}
\usepackage{geometry}
\usepackage{graphicx}
\usepackage{xcolor}
\usepackage{hyperref}

\newtheorem{theorem}{Theorem}[section]
\newtheorem{lemma}[theorem]{Lemma}

\theoremstyle{definition}
\newtheorem{definition}[theorem]{Definition}

\theoremstyle{remark}
\newtheorem{remark}[theorem]{Remark}
\DeclareMathOperator   {\tr}     {tr}

\title{\textbf{Singularities of two-dimensional Nijenhuis operators}}
\author{Dinmukhammed Akpan\\Email: dinmukhammed.akpan@uni-jena.de}
\date{}
\begin{document}
\maketitle

\section{Introduction}\label{sec1}

An important issue in differential geometry is the reduction of tensor fields to a constant form by a suitable coordinate transform. 

For example, the necessary and sufficient condition for the existence of coordinates in which a Riemannian metric has constant form (i.e., for the intergability of a metric) is the vanishing 
of the curvature tensor $R^i_{jkl}$, corresponding to this metric.

Similar condition for a differential 2-form $\omega$ is closedness ($(d\omega)_{ijk}=0$). 

The next natural object {  in the study of} tensor fields of rank two on smooth manifolds is an operator field $L$. 
The vanishing of the Nijenhuis { torsion}  
is a necessary condition {  for integrability (= reducibility to a constant form in some coordinates)} of an operator field. { It becomes also sufficient for semisimple operators with constant eigenvalues. However in general, both eigenvalues and algebraic type of a Nijenhuis operator may vary on $M$.} Nijenhuis operators also appear in various problems concerning integrable systems and other branches of mathematics. 

Many properties and applications of Nijenhuis operators were described in recent papers by A.\,V.~Bolsinov, V.\,S.~Matveev, A.\,Yu.~Konyaev \cite{bib1}, \cite{bib2}, \cite{bib3},\cite{bib4},\cite{bib5}. 
The necessary definitions are given in the next section.

The problem discussed in the paper was posed in the article by A.\,V.~Bolsinov, V.\,S.~Matveev, E.\,Miranda and S.\,L.~Tabachnikov \cite[Problem 5.13]{bib6} and {  is related to a description of singularities of Nijenhuis operators in dimension 2}. {  In elementary terms,} it sounds 
like this: describe all the functions $f (x, y)$ of two variables defined in a neighbourhood of
$(0,0) \in \mathbb{R}^2$ such that $f_y(0,0) = 0$ and $\frac{f_x(x - f_x) - f}{f_y}$ is smooth  { (more precisely, extends up to a smooth function at those points where $f_y=0$)}. {   The relation of this question to singularities of Nijenhuis operators in dimension 2 is based on the following fact:  if the trace and determinant of a Nijenhuis operator are functionally independent on $M^2$ almost everywhere, then $L$ can be uniquely reconstructed from them by means of a certain explicit formula \cite{bib1}.  This formula, however, involves division by the Jacobian of the characteristic map 
$$
\Phi: M^2 \to \mathbb R^2, \qquad  p \in M^2 \mapsto (\tr L(p), \det L(p))\in\mathbb R^2.
$$   
so that it does not work directly at singular points of $\Phi$ where $\operatorname{Jac}(\Phi)=0$  (in the context of Nijenhuis geometry, such points are known as {\it differentially degenerate}).  It means, in particular, that singularities of the characteristic map associated with a Nijenhuis operator must be very special.
If we assume in addition that $d \tr L \ne 0$ and take $x=\tr L$ as one of two local coordinates on $M^2$,  then the map $\Phi (x,y) = (x, f(x,y))$ may serve as a characteristic map of a smooth Nijenhuis operator if and only if $f(x,y)$ satisfies the above condition.  In other words,  answering the above problem is essentially equivalent to description of differentially degenerate singular points of two-dimensional Nijenhuis operators under additional condition that $d\tr L\ne 0$.  Recall that the structure of differentially non-degenerate singularities of Nijenhuis operators is well known, see \cite{bib1} and Theorem \ref{T1} below.}

The author thanks A.\,A.~Oshemkov, A.\,Yu.~Konyaev, A.\,V.~Bolsinov, V.\,S.~Matveev, E.\,A.~Kudryavtseva, and V.\,A.~Kibkalo for valuable advice and comments on the work.

\section{The main definitions}\label{sec2}

\begin{definition}[\cite{bib8}]

Let $M$ be a smooth manifold, $L$ a smooth operator field on it, and $v,\,u$ smooth vector fields. Then the Nijenhuis torsion $({  {\mathcal N_L}})^i_{jk}$ {  of $L$} is a tensor of type $(1,2)$ given by the formula
\begin{equation*}
{  {\mathcal N_L}}{  (v,u)} = [Lv, Lu] + L^2[v,u] - L[Lv,u] - L[v, Lu],
\end{equation*}
or in coordinates
\begin{equation*}
({  {\mathcal N_L}})^i_{jk} = L^l_j \frac{\partial L^i_k}{\partial x^l} - L^l_k \frac{\partial L^i_j}{\partial x^l} - L^i_l \frac{\partial L^l_k}{\partial x^j} + L^i_l \frac{\partial L^l_j}{\partial x^k}.
\end{equation*}
\end{definition}

\begin{definition}
{  An} operator $L$ is called a Nijenhuis operator if ${  {\mathcal N_L}}\equiv 0$.
\end{definition}

\begin{remark}

As \textit{smoothness} we will always assume infinite differentiability. 
\end{remark}

\begin{definition} [\cite{bib1}]
Let $M$ be an $n$-dimensional smooth manifold. A point $P\in M$ is called \textit{algebraically generic} with respect to an operator $L$ if there exists  a neighbourhood 
of this point in which the algebraic structure of $L$ does not change, that is, the structure of its Jordan normal form is preserved. Otherwise the point is called \textit{singular}.

A point $P\in M$ is called \textit{differentially non-degenerate} if the differentials of the coefficients of the characteristic polynomial of $L$ are linearly independent at this point. Otherwise the point is called \textit{differentially degenerate}.
\end{definition}

\begin{remark}
The values of the eigenvalues are not important for us, so in the {   real} two-dimensional case there are four possibilities:
\begin{equation*}
1. 
\begin{bmatrix}
\lambda_1 & 1 \\
0 & \lambda_1
\end{bmatrix}
\qquad
2. 
\begin{bmatrix}
\lambda_1 & 0 \\
0 & \lambda_2 
\end{bmatrix}
\qquad
{  3. 
\begin{bmatrix}
\lambda_1 & \lambda_2 \\
-\lambda_2 & \lambda_1
\end{bmatrix}
}
\qquad
4. 
\begin{bmatrix}
\lambda_1 & 0 \\
0 & \lambda_1
\end{bmatrix}
\end{equation*}

In case {  4}, $P$ is called of \textit{scalar type}.

\end{remark}

\begin{remark}
For the two-dimensional case, the Nijenhuis torsion has 8 components, but only $({  {\mathcal N_L}})^1_{1,2} = -({  {\mathcal N_L}})^1_{2,1}$ and $({  {\mathcal N_L}})^2_{1,2} = -({  {\mathcal N_L}})^2_{2,1}$ are essential, the rest are zero.
\end{remark}

\section{Singularities of Nijenhuis operators}\label{sec3}

From \cite{bib1} we know how to restore a Nijenhuis operator in neighbourhood of a differentially non-degenerate point from its invariants (coefficients of the characteristic polynomial).  In dimension 2,  the formula can be found in \cite{bib3}, \cite{bib6}.

\begin{theorem} [\cite{bib6}] \label{T1}
Let $L$ be a two-dimensional {  Nijenhuis} operator and $\tr L = x$, $\det L = f(x,y)$, with $f_y(x,y)\neq 0$. Then the point $(0,0)$ is differentially non-degenerate and $L$ is given in coordinates $(x,y)$ by the following formula:  
\begin{equation}
\label{eq:01}
L =
\begin{bmatrix}
x - f_x & -f_y \\
\frac{f_x(x - f_x) - f}{-f_y} & f_x
\end{bmatrix}.
\end{equation}
\end{theorem}

Notice that if in this theorem we choose the second coordinate to be $\tilde y=-f(x,y)=-\det L$, then  $L$ takes the standard canonical form  from \cite{bib1} for differentially non-degenerate points, namely $L=\begin{pmatrix} x & 1 \\ \tilde y & 0 \end{pmatrix}$.  

{  Formula \eqref{eq:01} holds at every differentially non-degenerate point.  Hence, if such points are everywhere dense (equivalently,  if $f_y\ne 0$ almost everywhere),  then $L$ can still be reconstructed from this formula by continuity.  However,   the function $f(x,y)$  must be very special in order for $L^2_1 =  \frac{f_x(x - f_x) - f}{-f_y}$ to be smooth.

We will answer the following question in dimension 2: what form can a smooth function $f(x,y)$ have if there exists a Nijenhuis operator $L$ in a neighborhood of the point $(0,0)$ such that $\tr L = x$ and $\det L = f(x,y)$ ?  }

The problem can be divided into 2 cases:

1. $\frac{\partial f}{\partial y} \equiv 0$, i.e., $f = f(x)$;

2. $\frac{\partial f}{\partial y}(0,0) = 0$ (but $\frac{\partial f}{\partial y}\not\equiv0$).

The case when $f_y(0,0) \neq 0$ {   is} completely described by {  T}heorem \ref{T1}.

Note again that the case when $\frac{\partial f}{\partial y}(0,0) = 0$ corresponds to differentially degenerate point $(0,0)$.

Thus, in this case{, if $\frac{\partial f}{\partial y} \not\equiv 0$,} the description of singularities of $L$ {  reduces} to studying the smoothness of the fraction $\frac{f_x(x - f_x) - f}{f_y}$.

Case 1 is different from the others, since here there is no question about smoothness of the fraction, due to the fact that $f_y \equiv 0$.

\begin{theorem}\label{T2}
Let $L$ be a two-dimensional Nijenhuis operator and $\tr L = x$, $\det L= f(x)$. Then either $ L = 
\begin{bmatrix}
x - \alpha & 0 \\
c(x,y) & \alpha
\end{bmatrix}
$
or  
$
L = 
\begin{bmatrix}
\frac{x}{2} & 0 \\
c(x,y) & \frac{x}{2}
\end{bmatrix}
$
, where $\alpha\in \mathbb{R}$, $c(x,y)$ is an arbitrary smooth function. 
In particular, 
$f(x) = \alpha x - \alpha^2$ or $f(x) =\frac{x^2}{4}.$
\end{theorem}

\begin{remark}
 {  Theorem \ref{T2} allows us to construct Nijenhuis operators with a rather complicated set $\mathsf {Sing}$ of singular points, since $\mathsf{Sing}$ consists of zeros of an arbitrary smooth function $c(x,y)\not\equiv 0$ for $L$ having the latter form (with diagonal elements $\frac x2$). On the contrary,   
if $\alpha \neq 0$, then $L= \begin{bmatrix}
x - \alpha & 0 \\
c(x,y) & \alpha
\end{bmatrix}$ is diagonalisable with distinct eigenvalues in a neighborhood of $(0,0)$ and this point is not singular. If $\alpha = 0$, then  $\mathsf {Sing}$
coincides with the line $\{ x=0\}$.
}

\end{remark}

\begin{proof}
Let  a Nijenhuis operator  have the form
\begin{equation*}
L = 
\begin{bmatrix}
a(x,y) & b(x,y) \\
c(x,y) & d(x,y)
\end{bmatrix}.
\end{equation*}
Then the essential components of the Nijenhuis torsion of $L$ are 
$$
({  {\mathcal N_L}})^1_{1,2} = a_y (a-d) + (bc)_y - (a+d)_x b, \quad
({  {\mathcal N_L}})^2_{1,2} = d_x (a-d) - (bc)_x + (a+d)_y c.
$$
Let us write the conditions of the theorem and write out the components of the Nijenhuis torsion:
\begin{equation*}
\begin{cases}
a + d = x \\
a_x + d_x = 1 \\
a_y + d_y = 0 \\
ad - bc = f(x) \\
a_y d +a d_y - b_y c - b c_y = 0 \\
({  {\mathcal N_L}})^2_{1,2} = f' - d = 0 \\
({  {\mathcal N_L}})^1_{1,2} = -b = 0
\end{cases}
\end{equation*}
From here we find explicit solutions for the components of the operator:
\begin{equation*}
\begin{cases}
a = x - f' \\
d = f' \\
b \equiv 0 \\
c - \text{arbitrary function} \\
(x - f')f' = f
\end{cases}
\end{equation*}
The last differential equation gives the form of the determinant:
\begin{equation*}
f(x) = \alpha x - \alpha^2 \qquad \text{or} \qquad f(x) =\frac{x^2}{4}{.}
\end{equation*}
{The latter form for $f(x)$ holds if the function $g(x)=\frac{\tr^2L}{4}-\det L=\frac{x^2}{4}-f(x)$ vanishes everywhere.
The former form for $f(x)$ holds in a neighbourhood of each point where the function $g(x)$ does not vanish. Since the function $g(x)$ is everywhere smooth and has the form $g(x)=(\frac x2-\alpha)^2$ in this neighbourhood, its zeros form the line $x=2\alpha$, thus the formulae for $g$ and $f$ hold everywhere.}
\end{proof}

Theorem \ref{T2} gives a complete answer to the first case of the problem.

\begin{remark}\label{RT2} 
In Theorem \ref{T2}, if $c(0,0)\ne 0$ 

then there exist new coordinates $(x,\tilde{y})$ centered at the origin, in which

$ L = 
\begin{bmatrix}
x - \alpha & 0 \\

1 & \alpha
\end{bmatrix}
$
, \ \
$
L = 
\begin{bmatrix}
\frac{x}{2} & 0 \\

1 & \frac{x}{2}
\end{bmatrix},
$
respectively.

\end{remark}

In the context of this paper, it will be more convenient to replace $f(x,y) = \det L$ with the function $g(x,y) =  \frac{\tr^2L}4-\det L$ which is nothing else but the discriminant of the characteristic polynomial of $L$.  

\begin{definition}
A function $g(x,y)$ is called \textit{admissible discriminant} if there exists a smooth {  Nijenhuis} operator $L$ such that $\tr L = x$ and $(\frac{\tr L}{2})^2 - \det L = g(x,y)$.
\end{definition}

In terms of the function $g(x,y)$, the general form of the Nijenhuis operator from Theorem \ref{T1} is rewritten as:
\begin{equation} \label {eq:1}
L =
\begin{bmatrix}
\frac{x}{2} + g_x & g_y \\
-\frac{g_x^2 - g}{g_y} & \frac{x}{2} - g_x
\end{bmatrix}.
\end{equation}

\begin{remark} \label{RAAF}
Due to Theorems \ref{T1} and \ref{T2},
$g(x,y)$ is an admissible discriminant  if and only if 
either $g(x,y)=(\frac x2-\alpha)^2$ for some $\alpha\in\mathbb R$, or $g(x,y)\equiv0$, or $g_y\not\equiv0$ and
$\frac{g_x^2 - g}{g_y}$ is smooth {  (more precisely, can be extended up to a smooth function at those points where $g_y =0$)}.
\end{remark}

The following theorem describes all Nijenhuis operators whose discriminant (as a smooth function) has a non-degenerate singularity at the point $(0,0)$.

\begin{theorem}\label{T3}
Let $g(x,y)$ be a smooth function such that 
$g(0,0) = g_y(0,0) ~=~ 0$, $g_{yy}(0,0)~\neq ~0$. Then $g(x,y)$ is an admissible discriminant if and only if  there exists a transformation of the second coordinate $\tilde y = \tilde y(x,y)$,  $\frac{\partial \tilde y}{\partial y}\ne 0$ that reduces $g$ to one of the following two forms: either $g=\pm\, {\tilde y}^2 + \frac{x^2}{4}$ or $g=\pm\, {\tilde y}^2$. In coordinates $(x, \tilde y)$, the corresponding Nijenhuis operators  are respectively
$
L_1^{\pm} ~=~
\begin{bmatrix}
x & \pm 2{\tilde y} \\
{
\frac{\tilde y}{2}} & 0
\end{bmatrix}
$
and
$
L_2^{\pm} ~=~
\begin{bmatrix}
\frac{x}{2} & \pm 2{\tilde y} \\
{
\frac{\tilde y}{2}} & \frac{x}{2}
\end{bmatrix}
$.

\end{theorem}

\begin{proof}
Since the function $g(0,y)$ has a non-degenerate singularity {at $y=0$}, it follows from {the} Morse parametric lemma that there exist {local} coordinates $\tilde x\equiv x$, $\tilde y=\tilde y(x,y)$ { centered at $(0,0)$} in which ${g} = \pm {\tilde y}^2 + \tau(x)${, where $\tau(x)$ is a smooth function with $\tau(0)=0$}. 
Smoothness condition of the function {$\frac{g^2_x - g}{g_y}$ means smoothness of the Nijenhuis operator (Remark \ref{RAAF}), which is independent on the choice of local coordinates preserving the coordinate $x$, thus it is equivalent to smoothness of the function}
\begin{equation*}
\frac{g^2_{\tilde x} - g}{g_{\tilde y}}=\frac{(\tau ')^2 - \tau  \mp {\tilde y}^2}{\pm 2{\tilde y}} = \frac{(\tau ')^2 - \tau}{\pm 2 {\tilde y}} - \frac{\tilde y}{2}{,}
\end{equation*}
{which} exactly means that $(\tau')^2 - \tau\equiv 0$. Hence $\tau(x) = (\frac{x}{2})^2$ or $\tau(x)\equiv 0$. Then, ${g} = \pm {\tilde y}^2 + \frac{x^2}{4}$ or ${g} = \pm {\tilde y}^2$ and we can recover Nijenhuis operator with formula~{\eqref{eq:1}}.
\end{proof}

\begin{figure}[h!]
\centering
\includegraphics[width=0.5\textwidth]{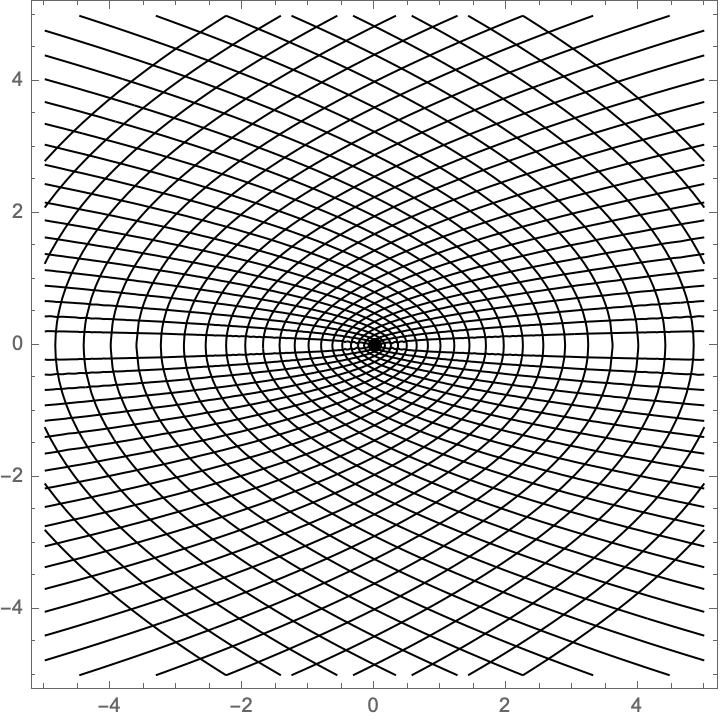}
\caption{The level lines of the real eigenvalues of the Niejenhus operator $L_1^{+}$ from Theorem~\ref{T3}.}\label{fig1}
\end{figure}

\begin{figure}[h!]
\centering
\includegraphics[width=0.5\textwidth]{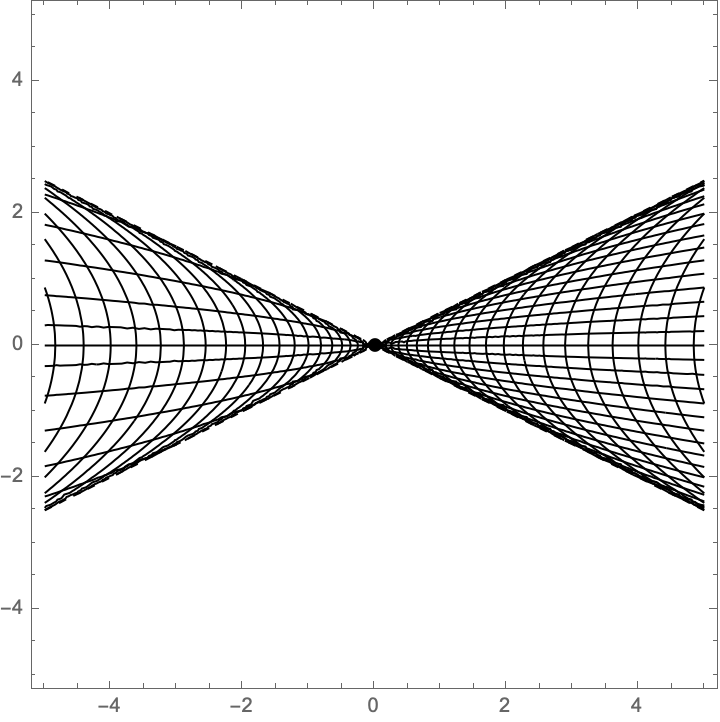}
\caption{The level lines of the real eigenvalues of the Niejenhus operator $L_1^{-}$ from Theorem~\ref{T3}. The domain with complex eigenvalues is given by the inequality $x^2 - 4y^2 < 0$.}\label{fig2}
\end{figure}

\begin{figure}[h!]
\centering
\includegraphics[width=0.5\textwidth]{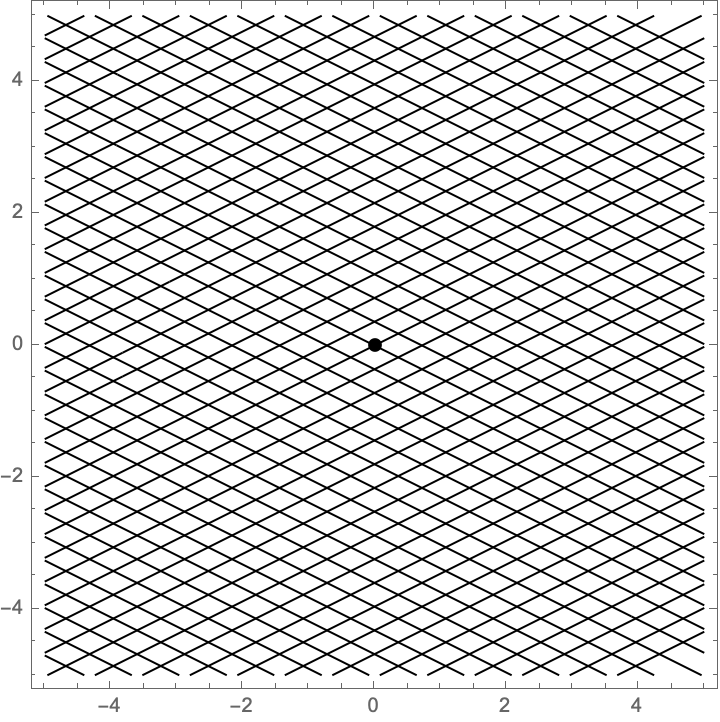}
\caption{The level lines of the real eigenvalues of the Niejenhus operator $L_2^{+}$ from Theorem~\ref{T3}. For $L_2^{-}$ all eigenvalues are complex.}\label{fig3}
\end{figure}

\begin{remark}{ 
We emphasize once again that Theorem {\color{red}\ref{T3}} gives a complete description of the Nijenhuis operators such that $\tr L =x$ and  $\det L$  (or equivalently, the discriminant of $L$)  has a Morse singularity in the second variable $y$.   In other words, it describes Nijenhuis operators  whose characteristic map $(\tr L, \det L): M \to \mathbb R^2$ has a non-degenerate singularity of rank one (called a fold point \cite{bib9}). In particular, this theorem shows that such operators are linearizable in the sense of the paper \cite{bib4}, where the linearization problem was discussed in the context of Nijenhuis geometry.}
\end{remark}

The next natural step is to consider Nijenhuis operators whose admissible discriminant has a degenerate singularity in the second variable. The answer to this question is given by Theorem 4, where the admissible discriminant has a cubic singularity in the second variable.

\begin{theorem}{\label{T4}}
Let $g(x,y)$ be an analytic function such that $g(0,0) = g_y(0,0) =g_{yy}(0,0) = 0$, $g_{yyy}(0,0)\neq 0$. Then $g(x,y)$ is an admissible discriminant if and only if { 
there exists a transformation of the second coordinate $\tilde y = \tilde y(x,y)$,  $\frac{\partial \tilde y}{\partial y}\ne 0$ that reduces $g$ to one of the following two forms: either
${g} = \tilde y^3$  or $g = \tilde y^3 + \frac{x^2}{4}$. In coordinates $(x, \tilde y)$, the corresponding Nijenhuis operators  are respectively: $
L_1 ~=~
\begin{bmatrix}
\frac{x}{2} &  3{\tilde y}^2 \\
\frac{\tilde y}{3} & \frac{x}{2}
\end{bmatrix}
$
and
$
L_2 ~=~ 
\begin{bmatrix}
x & 3{\tilde y}^2 \\
\frac{\tilde y}{3} & 0
\end{bmatrix}
$.
}
\end{theorem}

\begin{figure}[h!]
\centering
\includegraphics[width=0.5\textwidth]{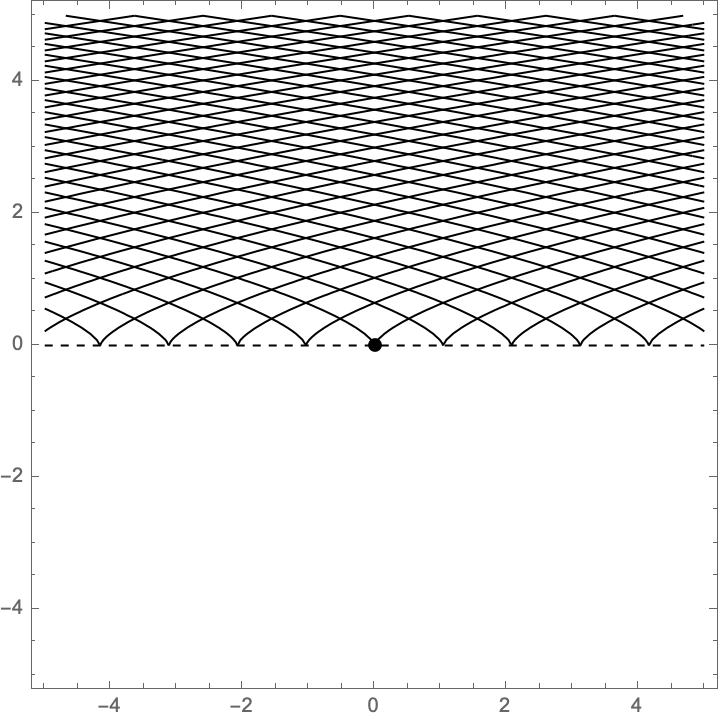}
\caption{The level lines of the real eigenvalues of the Niejenhus operator $L_1$ from Theorem~\ref{T4}.The domain with complex eigenvalues is given by the inequality $y < 0$.}\label{fig4}
\end{figure}

\begin{figure}[h!]
\centering
\includegraphics[width=0.5\textwidth]{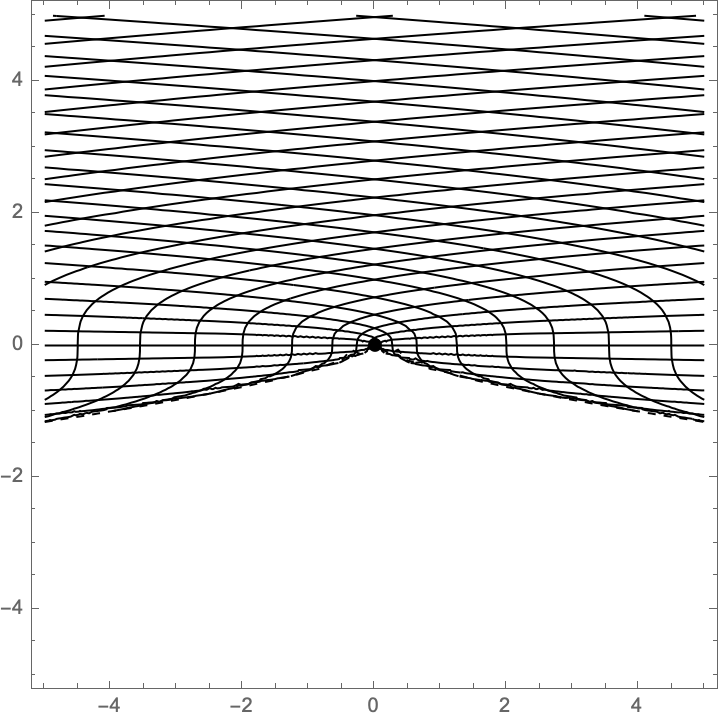}
\caption{The level lines of the real eigenvalues of the Niejenhus operator $L_2$ from Theorem~\ref{T4}. The domain with complex eigenvalues is given by the inequality $x^2 + 16y^3 < 0$.}\label{fig5}
\end{figure}

\begin{lemma}\label{L1}
The following formulas are true
\begin{equation*}
\frac{\partial^{(2k)}}{\partial y^{2k}} \Bigl (\frac{1}{y^2 + \tau(x)}\Bigr)\mid_{y=0} = \frac{(-1)^k(2k)!}{\tau^{k+1}(x)}{,}
\end{equation*}
\begin{equation*}
\frac{\partial^{(2k + 1)}}{\partial y^{2k+1}}\Bigl(\frac{y}{y^2 + \tau(x)}\Bigr)\mid_{y = 0} =  \frac{(-1)^{k}(2k + 1)!}{\tau^{k + 1}(x)}.
\end{equation*}
\end{lemma}

\begin{proof}
To prove the first formula, decompose the function into the series: $\frac{1}{y^2+\tau(x)} = \frac{1}{\tau(x)} - \frac{y^2}{\tau^2(x)} + \frac{y^4}{\tau^{3}(x)} { - } \ldots$. 

Then we see that the coefficients at even degrees of $y$ have the desired form.
The proof of the second formula is obtained similarly after multiplying the function from the first formula by~$y$.
\end{proof}

\begin{proof}
From the theorem on the general form of smooth functions in a neighborhood of zero with a degenerate singularity (see \cite{bib7}), 
we obtain that the function $g$ can be reduced by replacing the coordinates {with $\tilde x\equiv x$, $\tilde y=\tilde y(x,y)$} to the form ${g} = {\tilde y}^3 + \tau({\tilde x}) {\tilde y} + \beta({\tilde x})$ {where $\tau(0)=\beta(0)=0$}. 
Let's check the smoothness condition for the function
\begin{equation*}
\frac{g_{\tilde x}^2 - g}{g_{\tilde y}} = -\frac{\tilde y}{3} + \frac{ (\tau')^2}{3} + \underbrace{\frac{1}{3} \frac{{\tilde y}(6 \tau' \beta' - 2\tau) - \tau (\tau')^2 + 3(\beta')^2 - 3\beta}{3{\tilde y}^2 + \tau({\tilde x})}}_{I}.
\end{equation*}

It is clear that the question of the smoothness of the fraction $\frac{g_{\tilde x}^2 - g}{g_{\tilde y}}$ is equivalent to the study of the smoothness of $I$. 
We will look at the derivatives of this fraction with respect to $\tilde y$ at $\tilde y = 0$. We have an expression of the form $\frac{{\tilde y}A(x) + B(x)}{3{\tilde y}^2 +\tau(x)}$. 
Since the functions $A(x)$ and $B(x)$ are obtained from $\tau(x)$ and $\beta(x)$ by using arithmetic operations and differentiation, then $A(x)$ and $B(x)$ have a finite order at zero{, provided that they are not identical $0$}.
Denote by $m$, $n$ the orders of zero for $A(x)$ and $B(x)$ respectively, also denote by $l$ the order of zero for $\tau(x)$ {(we remark that if $A(x)\not\equiv0$ and $B(x)\not\equiv0$ then $\tau(x)\not\equiv0$, since $I$ is smooth)}.
Differentiating the fraction $I$ by $\tilde y$ and then assuming ${\tilde y} = 0$ we obtain from Lemma~\ref{L1} that the power of $\tau(x)$ in the denominator increases. 
It is clear that there is a number $N$ such that $Nl > m$ and $Nl > n$. Therefore, if the numerator of the fraction $I$ is nonzero, then its derivative of some order with respect to $\tilde y$ 
will not be continuous, that is, the only possibility is $A(x)\equiv B(x)\equiv 0$.

The conditions $A(x)\equiv B(x)\equiv 0$ in our case give a system of nonlinear ODE's{:}
\begin{equation*}
\begin{cases}
3 \tau' \beta' - \tau = 0 {,} \\
- \tau (\tau')^2 + 3(\beta')^2 - 3\beta = 0{.}
\end{cases}
\end{equation*}

Solving the system (by evaluating $\beta'$ in terms of $\tau,\tau'$ from the first equation, substituting it into the second equation and evaluating $\beta$ in terms of $\tau,\tau'$), {we get three solutions} 
$(\tau(x)\equiv 0, \beta(x) = \frac{x^2}{4})$, $(\tau(x) = -\frac{3x^{\frac{4}{3}}}{2^{\frac{8}{3}}}, \beta(x) = \frac{x^2}{8})$, and 
$(\tau(x) \equiv 0, \beta(x)\equiv 0)$. The second solution is not smooth, therefore it does not {satisfy} the conditions of the theorem and we can reestablish Nijenhuis operator with formula~{\eqref{eq:1}}.
Thus, we have found all possible solutions for third order singularities. 
\end{proof}

\begin{remark}
Figures 1-4 illustrate the level lines of real eigenvalues of the Nijenhuis operators described in Theorems 3 and 4.
\end{remark}

\section{Examples of Nijenhuis operators}\label{sec4}

\begin{theorem}\label{T5}
Let $g(x,y) = G(y)$ be a function that has zero of finite order $k \geq 2$. Then $g(x,y) = G(y)$ is an admissible discriminant.  {  The corresponding Nijenhuis operator is $L=\begin{pmatrix}  \frac{x}{2}\!\! \phantom{\Bigl(} & G_y \\ \frac{G}{G_y} & \frac{x}{2}  \end{pmatrix}$.}  
\end{theorem}

\begin{proof}
It is enough to check the smoothness condition for $\frac{g_x^2 - g}{g_y}= -\frac{G}{G_y}$. By {the} Hadamard lemma we have $G(y) = y^k F(y)$, where $k \geq 2$, $F(0)\neq 0$. 

Therefore, $\frac{G}{G_y} = y \frac{F(y)}{kF(y)+yF'(y)} \in C^{\infty}${, which allows us to define a Nijenhuis operator by the formula~\eqref{eq:1}}.
\end{proof}

\begin{remark}
The condition of finiteness of order of zero in Theorem \ref{T5} is not superfluous. 
Indeed, for a flat function $f$ (i.e., such that $f(0) = f'(0) = \ldots = f^{(k)}(0) {=} \ldots = 0$), the situation is more complicated since it is possible to construct an example 
of a flat function for which the smoothness condition for $\frac{g}{g_y}$ does not hold. 
For example, if $g(y) = e^{-\frac{1}{\sqrt{y}}}$, then

$\frac{g}{g_y} =  2 \sqrt{y} y$ is not smooth at zero.
\end{remark}

\begin{remark}
In theorem \ref{T5} we can change {  the $y$-coordinate to reduce $G$ to the form} $G(y) = {\pm} y^k$, where $k$ is the order of zero. Then, in this case the Nijenhuis operator will be: $L = \begin{bmatrix} \frac{x}{2} & {\pm} k y^{k-1} \\ \frac{y}{k} & \frac{x}{2} \end{bmatrix}$. 
\end{remark}

\begin{theorem} \label{T6}
Let $g(x,y) = a(x) + b(x)y$, where $a(x)$ and $b(x)$ are smooth and {  have zeros of order $m$ and $k$ respectively at $x=0$}. Suppose $m\geq k\geq 2$. Then $g(x,y)$ is an admissible discriminant and the corresponding Nijenhuis operator is: 
$$L = \begin{bmatrix}
\frac{x}{2} + (b(x)C(x,y))_x & b(x)C_{y}(x, y) \\ 
-\frac{((b(x)C(x, y))_x )^2 - b(x)C(x, y)}{b(x) C_{y}(x, y)} & \frac{x}{2} - (b(x)C(x, y))_x 
\end{bmatrix},$$ where $C(x, y) = \frac{a(x)}{b(x)} + y$.
\end{theorem}

\begin{remark}
In Theorem \ref{T6} we can take $\frac{a(x)}{b(x)} + y$ as a new coordinate $\tilde{y}$. Then 
 $L = \begin{bmatrix}
\frac{x}{2} + b'(x) \tilde{y} & b(x) \\ 
-\frac{(b'(x) \tilde y)^2 - b(x) \tilde y}{b(x)} & \frac{x}{2} - \tilde y
\end{bmatrix}$ in coordinates $x,\tilde y$.
\end{remark}

\begin{proof}
As we already know, {  the admissibility of $g(x,y)$ is equivalent to the smoothness of the fraction $\frac{g_x^2 - g}{g_y}$}. 

In order to check that the fraction 
\begin{equation*}
\frac{(b')^2y^2 + y(2a'b' - b) + ((a')^2 - a)}{b}
\end{equation*}
is smooth, it is enough to show that all fractions $\frac{(b')^2}{b}$, $\frac{{a}'b'}{b}$, $\frac{b}{b}$, $\frac{(a')^2}{b}$, $\frac{a}{b}$ are smooth, but this follows from the condition on orders of zero.
\end{proof}

One of the interesting classes of functions is the space of homogeneous polynomials.

\begin{theorem}\label{T7}
Let  $g(x,y)$ be a homogeneous polynomial {of degree $\ne2$, with $g_y\not\equiv0$}, then $g(x,y)$ is an admissible discriminanti n a neighbourhood of the origin if and only if  $g(x,y) = {\pm} x^m(ax +by)^k$, where $b\neq 0$, $k\geq {1}$, $m\geq 0${, $m\ne1$}.

\end{theorem}

\begin{remark}
{The condition of degree $\ne2$ in Theorem \ref{T7} is not superfluous, as can be seen from the example in Theorem \ref{T3}.}
Theorem \ref{T7} can be reformulated as follows: let $g(x,y)$ be a homogeneous polynomial {of degree $\ne2$} with $g_y\not\equiv 0$, then $g(x,y)$ is an admissible discriminant if and only if ${g} = {\pm} x^m\tilde{y}^k$ for some linear combination $\tilde y = ax + by$, $b\ne 0$. The corresponding Nijenhuis operator in coordinates $x,\tilde y$ will be: $L ~=~ \begin{bmatrix}
\frac{x}{2} {\pm} mx^{m-1} {\tilde y}^k & {\pm} k {x^m} {\tilde y}^{k-1} \\
\frac{{\pm m^2} x^{m-2} {\tilde y}^{k+1} - {\tilde y}}{k} & \frac{x}{2} {\mp} mx^{m-1} {\tilde y}^k
\end{bmatrix} $ where $k\ge{1}$, $m\ge 0${, $m\ne1$}.

\end{remark}

\begin{proof}
Let $m\geq 2$, check the smoothness of $\frac{g_x^2 - g}{g_y}$ with the explicit form $g(x,y)${:}
\begin{equation*}
g_x = {\pm} mx^{m-1}(ax + by)^k {\pm} kax^m(ax + by)^{k-1} \
\end{equation*}
\begin{equation*}
\Rightarrow \ (g_x)^2 = x^{2m - 2}(ax + by)^{2k - 2}(m(ax + by) + kax)^2{.}
\end{equation*}
Thus
\begin{equation*}
\frac{(g_x)^2 - g}{g_y} = \frac{ x^{2m - 2}(ax + by)^{2k - 2}(m(ax + by) + kax)^2 {\mp} x^m(ax + by)^k }{{\pm}kbx^m(ax + by)^{k-1}}
\end{equation*}

\begin{equation*}
= \frac{{\pm} x^{m-2} (ax + by)^{k-1} [m(ax + by) + kax]^2 - (ax + by)}{kb} 
{\in C^\infty.}
\end{equation*}

Let the fraction $\frac{g_x^2 - g}{g_y}$ be smooth, and let $g(x,y)$ be a homogeneous polynomial. Let's fix the degree of this polynomial $\deg(g(x,y)) = {n\ne2}$, then $\deg(g_x) = \deg(g_y) = {n} -1$. 
The fact that $\frac{g_x^2 - g}{g_y}$ is smooth exactly means that it is a polynomial, and:
\begin{equation*}
\frac{(g_x)^2 - g}{g_y} = \text{polynomial of degree ({$n - 1$})} \ + \ \text{linear polynomial}{.}
\end{equation*}
Now consider $g(x,y)$ which is a polynomial in $y$ with coefficients from $x$. Then the fact that $\frac{g}{g_y}$ is a linear function implies that $GCD(g, g_y) = g_y$. 
By Bezout's theorem, we have an explicit form of such polynomials{:}
\begin{equation*}
g(x,y) = {\pm} x^m(ax+by)^k{.}
\end{equation*}
Let $m = 0$, then $g(x,y) = {\pm} (ax +by)^k${,} hence $\frac{g_x^2 - g}{g_y}=\frac{{\pm a^2k^2(ax+by)^{k-1}} -(ax+by)}{bk}{\in C^\infty}.$
Let $m = 1$, then $g(x,y) = {\pm} x(ax + by)^k${,} hence $\frac{g_x^2 - g}{g_y}={\frac{{\pm} x^{-1} (ax + by)^{k-1} [(ax + by) + kax]^2 - (ax + by)}{kb}}\not\in C^\infty.$

\end{proof}

\bibliography{sn-bibliography}

\end{document}